\documentclass[12pt]{article}
\usepackage{amssymb}
\usepackage{amsthm}
\usepackage{graphicx}
\usepackage{psfrag}
\usepackage{graphics}
\usepackage{amsthm}

\input pictex.tex
\newcommand{\clo}{\mathrm{S}^1}

\theoremstyle{definition}

\newtheorem{thm}{Theorem}[section]

\newtheorem{prop}[thm]{Proposition}

\newtheorem{lem}[thm]{Lemma}

\newtheorem{rem}[thm]{Remark}

\input epsf

\date{}
\author{Andr\'es Navas}

\begin{document}

\title {Reduction of cocycles and groups of diffeomorphisms of the circle}
\maketitle

\vspace{-0.5cm}

\noindent{\bf Abstract.} We prove two theorems of reduction of cocycles taking values in 
the group of diffeomorphisms of the circle. They generalise previous results obtained by 
the author concerning rigidity for smooth actions on the circle of Kazhdan's groups and higher 
rank lattices.\\

\vspace{0.05cm}

\noindent{\bf Subject classification AMS (2000)}: primary 57S20; secondary 37A20.\\

\vspace{0.05cm}


%%%%%%%%%%%%%%%%%%%%%%%%%%%%%%%%%%%%%%%%%%%%%%%%%%%%%%%%%%%%%%%%%%%%%%%%%%%%%%%%%%%%%%%%%%%

\section*{Introduction}

\hspace{0.45cm} As a general principle, it should be possible to extend all results involving 
groups satisfying Kazhdan's property (T) to the setting of {\em cocycles} introduced by 
R. Zimmer (see \S 1 for definitions). We list below several examples. 
In all of them, $G$ is a Kazhdan group acting ergodically by measure 
preserving transformations on a probability space $\Omega$, and 
$\alpha : G \times \Omega \rightarrow H$ is a Borel cocycle 
taking values in a locally compact topological group $H$. 

\noindent{(a) If the group $H$ is amenable, then $\alpha$ is cohomologous to a cocycle taking 
values into a compact subgroup of $H$ \cite{Zimmer-libro}. This result has been extended in 
\cite{jol} to the case when $H$ satisfies Haagerup's property, {\em i.e.} when it is a-(T)-menable 
(see also \cite{An}).}

\noindent{(b) If $H$ is a Lie group, then $\alpha$ is cohomologous to a cocycle taking 
values into a Kazhdan subgroup of $H$ \cite{Zim-inv}.}

\noindent{(c) If $H$ is the group of isometries of a real tree, then $\alpha$ is cohomologous 
to a cocycle taking values into the stabiliser of some point of the tree \cite{AS}. This 
result can be easily extended to the case in which $H$ is the group of automorphisms 
of a measured wall-space \cite{CJV}.}

\vspace{0.15cm}

In this article we add to the list above the following result.

\vspace{0.35cm}

\noindent{\bf Theorem A.} {\em Let $\alpha : G \times \Omega \rightarrow 
\mathrm{Diff}_+^{1+\tau}(\clo)$ be a Borel cocycle, 
where $\tau > 1/2$ and $G$ is a compactly generated 
topological group whose action on $\Omega$ is measure preserving and 
ergodic. Suppose that, for each $g \in G$, the map $x \mapsto \alpha(g,x)$ takes values 
a.e. into a bounded subset of $\mathrm{Diff}_+^{1+\tau}(\clo)$. If $G$ has Kazhdan's 
property (T), then as a cocycle into $\mathrm{Homeo}_+(\clo)$, $\alpha$ is cohomologous 
to a cocycle taking values into the group of (euclidean) rotations.}

\vspace{0.35cm}

Recall that, for $\tau > 1/2$, every subgroup of $\mathrm{Diff}_+^{1+\tau}(\clo)$ that 
satisfies Kazhdan's property (T) is conjugated to a (compact) group of Euclidean 
rotations. This result was obtained by the author in \cite{Na-kazh} (see also \cite{Rez}). 
Theorem A above is just a generalisation of this fact into the framework of cocycles. It 
should be mentioned that a better result is known when $G$ is a Lie group satisfying 
Kazhdan's property (T). Indeed, using ideas introduced by \'E. Ghys \cite{ghys}, 
D. Morris Witte and \hspace{20cm} R. Zimmer proved  
\cite{WZ} a similar result in that case for cocycles into $\mathrm{Diff}^1_+(\clo)$ 
(and also a partial result for cocycles into $\mathrm{Homeo}_+(\clo)$). 

The reader which is familiar to ideas and techniques from rigidity theory for 
lattices in semisimple Lie groups could think that the proof of Theorem A is 
just a translation of the proof given in \cite{Na-kazh} to the language of
cocycles. Nevertheless, there are several (more than technical) problems which appear 
in that translation, and this is the main reason that motivated the author to write 
this article. Similar problems appear when we deal with the higher rank case. For that case 
we follow essentially the approach of \cite{Shalom} (see also \cite{Monod-Shalom}): 
our prototype of higher rank group will be a product $G = G_1 \times \cdots \times G_k$ 
of $k\geq 2$ compactly generated topological groups $G_i$. We will also suppose that 
$G$ acts on $\Omega$ {\em ergodically irreducibly}, that is the action of each subgroup 
$G_i' = G_1 \times \cdots \hat{G_i} \times \cdots \times G_k$ on $\Omega$ is ergodic. 
A typical example of this situation is given by the action of $G$ on the quotient 
$G / \Gamma$, where $\Gamma$ is a lattice in $G$ whose projections
into each $G_i'$ are dense. The following result is inspired by \cite{Na-superrig}. The 
precise meaning of its statement will be clarified in \S 3.

\vspace{0.35cm}

\noindent{\bf Theorem B.} {\em Let $\alpha : G \times \Omega \rightarrow 
\mathrm{Diff}_+^{1+\tau}(\clo)$ be a Borel cocycle, where $G=G_1 \times \cdots \times G_k$ 
is a product of $k \geq 2$ compactly generated topological 
groups and $\tau > 1/2$. Suppose that the $G$-action on 
$\Omega$ is ergodically irreducicle, and that for each $g \in G$ the map $x \mapsto \alpha(g,x)$ 
takes values a.e. into a bounded subset of $\mathrm{Diff}_+^{1+\tau}(\clo)$. If the action of 
$G$ on $\Omega \times \clo$ does not preserve any probability measure, then up to a topological 
semiconjugacy and a finite cover, $\alpha$ is cohomologous in 
$\mathrm{Homeo}_+(\mathrm{S}^1)$ to a cocycle given by a homomorphism 
from $G$ into the group of direct homeomorphisms of the circle which factors through one of 
the $G_i$. Moreover, if each $G_i$ is non discrete and almost topologically simple, then the 
image of this homomorphism coincides with some finite cover of $\mathrm{PSL}(2,\mathbb{R})$.}

\vspace{0.35cm}

It must be mentioned that Theorem B has a version in class $C^1$ (and 
also a partial version 
in class $C^0$) when $G$ is a higher rank semisimple Lie group without 
Kazhdan's property (T). This is stated explicitely in \cite{WZ}, and the proof is 
based on some of the arguments of \cite{ghys}.

The plan of this paper is the following. In \S 1 we recall some definitions and the main technical 
tools used in the proof of theorems A and B, namely the Cocycle
Reduction Lemma (due to R. Zimmer) and the Superrigidity Theorem 
for Reduced Cohomology (due to Y. Shalom) respectively. 
In \S 2 we extend the technique introduced in \cite{Na-kazh} to the 
case of cocycles. The main step for the proof of Theorem A is given in \S \ref{liso}, where 
we prove that the action of the group of orientation preserving homeomorphisms of the circle 
on some space of ``stable" geodesic currents is {\em smooth} (that is, the orbits are locally 
closed). This is a fundamental fact that must be verified in order to apply cocycle reduction 
and to finish the proof of Theorem A. Finally, the proof of Theorem B is given in \S 3. 

\vspace{0.2cm}

\noindent{\bf Remark.} Theorems A and B (and also the results from \cite{Na-superrig} and 
\cite{Na-kazh}) are still true for cocycles taking values in the group
of diffeomorphisms of $\clo$ of slightly lower differentiability class 
than $C^{3/2}$. However, we will discuss this point only in the 
Appendix, since the proofs rely on recent developments on property (T)
and the geometry of ``almost'' Hilbert spaces, and are independent of the rest of the paper. 

%%%%%%%%%%%%%%%%%%%%%%%%%%%%%%%%%%%%%%%%%%%%%%%%%%%%%%%%%%%%%%%%%%%%%%%%%%%%%%%%%%%%%%%%%%%%%%%%

\section{Some preliminary facts}

%%%%%%%%%%%%%%%%%%%%%%%%%%%%%%%%%%%%%%%%%%%%%%%%%%%%%%%%%%%%%%%%%%%%%%%%%%%%%%%%%%%%%%%%%%%%%%%%

\subsection{Zimmer's Cocycle Reduction Lemma} 

\hspace{0.45cm} Let $\Omega$ be a Borel space endowed with a probability measure $\mu$ 
and let $G$ be a group acting on $\Omega$ and preserving the measure class of $\mu$. A 
Borel map $\alpha : G \times \Omega \rightarrow H$ taking values into a topological 
group $H$ is called a {\em cocycle} if for all $g_1,g_2 \in G$ and a.e. $x \in \Omega$, 
$$\alpha(g_1g_2,x)=\alpha(g_1,g_2(x)) \alpha(g_2,x).$$
When $H$ is identified to a group of automorphisms of some space $\Omega'$, the fact that 
a (Borel) map $\alpha : G \times \Omega \rightarrow H$ is a cocycle is equivalent to 
that the map from $G$ into the group of automorphisms of $\Omega \times \Omega'$ given 
by $g(x,y) = (g(x),\alpha(g,x)(y))$ is a homomorphism. 

Two cocycles $\alpha$ and $\beta$ are {\em cohomologous} if there exists a Borel map 
$\varphi : \Omega \rightarrow H$ such that, for a.e. $x \in \Omega$, 
$$\alpha(g,x) = \varphi(g(x))^{-1} \beta (g,x) \varphi(x).$$
When $H$ is the group of automorphisms of a space $\Omega'$, a Borel function 
$\psi : \Omega \rightarrow \Omega'$ is said to be {\em equivariant} if for all $g \in G$ 
one has $\psi(g(x)) = \alpha(g,x) \psi(x)$ for a.e. $x \in \Omega$. In what follows 
we will always suppose that the Borel structure of $H$ is countably generated 
(\cite{Zimmer-libro}, page 10).

\vspace{0.25cm}

\noindent{\bf Cocycle Reduction Lemma (\cite{Zimmer-libro}, page 108).} {\em Let 
$\alpha : G \times \Omega \rightarrow H$ be a Borel cocycle, $G$ acting ergodically on 
$\Omega$. Suppose that $\Omega'$ is a continuous $H$-space on which the action is 
smooth (i.e. the orbits are locally closed). If there exists an equivariant function 
from $\Omega$ to $\Omega'$, then there is a point $y \in \Omega'$ such that $\alpha$ 
is cohomologous to a cocycle taking values into the stabiliser $H_y$ of $y$ in $H$.}

\vspace{0.25cm}

The proof of this lemma is very simple: it is based on the principle that functions which are 
constant along the orbits of an ergodic action are essentially constant. In order to apply 
this principle, it is important to have a good structure for the space of orbits by $H$ in 
$\Omega'$. Actually, R. Zimmer states the lemma only for locally compact groups 
$H$, but it is easy to see that this hypothesis is not completely
necessary. The essential hypothesis are that the Borel structure 
of $H$ is countably generated and that the $H$-action on $\Omega'$ is smooth. 
We clarify all of this because we will deal in the sequel with groups
of homeomorphisms of the circle, which in general are not locally compact \cite{ghys2}.

%%%%%%%%%%%%%%%%%%%%%%%%%%%%%%%%%%%%%%%%%%%%%%%%%%%%%%%%%%%%%%%%%%%%%%%%%%%%%%%%%%%%%%%%%%%%%%%%

\subsection{Kazhdan's property (T)} 

\hspace{0.45cm} Recall that a locally compact topological group $G$ has Kazhdan's property (T) 
(or to simplify, is a Kazhdan group) if every (continuous) action of $G$ by isometries of a 
separable Hilbert space has an invariant vector. Actually, this is not the original 
definition given by D. Kazhdan, but an equivalent one due to 
J. P. Serre and denoted property (FH). In this direction, one can naturally 
define property (FH) for actions, and Theorem A will still be true in this more general 
setting when $G$ fails to be a Kazhdan's group but its action on
$\Omega \times \clo$ has property (FH). 
Let us remark that properties (T) and (FH) are also equivalent for 
actions, and more generally for equivalence relations \cite{An}. 

%%%%%%%%%%%%%%%%%%%%%%%%%%%%%%%%%%%%%%%%%%%%%%%%%%%%%%%%%%%%%%%%%%%%%%%%%%%%%%%%%%%%%%%%%%%%%%%%

\subsection{Shalom's Superrigidity Theorem} 

\hspace{0.45cm} A continuous action of a locally compact topological group $G$ by isometries of 
a metric space $(X,d)$ is called {\em uniform} if there exist $\varepsilon > 0$ and a compact 
generating set $\mathrm{C}$ of $G$ such that for all $x \in X$ there exists $g\in\mathrm{C}$ 
satisfying $d(g(x),x) \geq \varepsilon$. We will be mainly interested in actions by (affine) 
isometries of a (real) Hilbert space $\mathcal{H}$. Recall that the group of isometries of 
such a space is the semidirect product between the unitary group and the group of 
translations. For an isometry $A$ we will denote by $\theta \in U(\mathcal{H})$ and 
$c \in \mathcal{H}$ its unitary and translation component respectively. 

Among the remarkable results obtained by Y. Shalom in \cite{Shalom}, 
for future reference we state here the Superrigidity Theorem for Reduced 
Cohomology, which will be essential in \S 3. 

\vspace{0.2cm}

\noindent{\bf Theorem \cite{Shalom}.} Let $G = G_1 \times \cdots \times G_k$ be a topological 
group which is the product of $k \geq 2$ compactly generated groups $G_i$. If $A = \theta + c$ 
is a uniform isometric action of $G$ on a Hilbert space $\mathcal{H}$, then there exists a 
non zero vector in $\mathcal{H}$ which is fixed by one of the $\theta(G_i')$, where 
$G_i' = G_1 \times \cdots \times \hat{G_i} \times \cdots \times G_k$.  

\vspace{0.2cm}

We will also use the following classical lemma, due to P. Delorme \cite{HV}.

\vspace{0.2cm}

\noindent{\bf Lemma.} {\em Let $A = \theta + c$ be an action of a locally compact topological 
group $G$ by isometries of a Hilbert space $\mathcal{H}$. If $A$ has no global fixed point 
and is non uniform, then there exists a sequence $(K_n)$ of unitary vectors in $\mathcal{H}$ 
which is $\theta$-almost invariant, that is for each compact subset $\mathrm{C}$ of $G$ the 
value of \hspace{0.05cm} $\sup_{g \in \mathrm{C}} \| \theta(g)K_n - K_n \|$ \hspace{0.03cm} 
tends to zero as $n$ tends to infinity.}

%%%%%%%%%%%%%%%%%%%%%%%%%%%%%%%%%%%%%%%%%%%%%%%%%%%%%%%%%%%%%%%%%%%%%%%%%%%%%%%%%%%%%%%%%%%%%%%

\subsection{Stable geodesic currents and convergence groups}
\label{conv}

\hspace{0.45cm} Recall that a geodesic current on the circle is a Radon measure defined on 
the space $\clo \times \clo \setminus \Delta$ which is invariant by the flip $(u,v) \mapsto (v,u)$. 
($\Delta$ denotes the corresponding diagonal.) We say that such a current $\nu$ is 
{\em stable} if 
$$\nu([a,b] \times [c,c]) = 0 \qquad \mbox { and } \qquad 
\nu([a,b[ \times ]b,c]) = \infty \quad \mbox{ for all } \quad a\!<\!b\!<\!c\!<\!a.$$
The notion of stable geodesic current was introduced in \cite{Na-kazh}, where the following 
proposition was stated without proof. The author is indebted to
J. C. Yoccoz for the appointment of the argument below, which simplifies 
an original one given in a previous version of this paper.

\vspace{0.2cm}

\begin{prop} {\em If $\nu$ is a stable geodesic current, then the group $G_{\nu}$ 
of (orientation preserving) homeomorphisms of the circle whose diagonal action on 
$\clo \times \clo$ preserves $\nu$ is topologically conjugated to a subgroup of 
$\mathrm{PSL}(2,\mathbb{R})$.}
\label{conver}
\end{prop}

\noindent{\bf Proof.} We will prove that $G_{\nu}$ has the Convergence 
Property, that is, each sequence $(g_n)$ in $G_{\nu}$ is either equicontinuous 
or contains a subsequence $(g_{n_k})$ such that for some points $a, b \in \clo$ 
one has that $g_{n_k}(x)$ tends to $b$ for all $x \in \clo \setminus \{a\}$ 
and that $g_{n_k}^{-1}(x)$ tends to $a$ for all $x \in \clo \setminus \{b\}$. 
By \cite{CJ,Ga,Hi3} and \cite{Tu}, this implies that $G_{\nu}$ is 
topologically conjugated to a subgroup of $\mathrm{PSL}(2,\mathbb{R})$.

Let us suppose that $(g_n)$ is a non equicontinuous sequence in $G_{\nu}$. Up to a 
subsequence, we may assume that there exists a point $a \in \clo$ and two sequences 
$(x_n)$ and $(y_n)$ converging to this point, such that $x_n < a < y_n < x_n$ for 
all $n \in \mathbb{N}$, and such that $g_n(x_n)$ tends to some point 
$x_{\infty} \in \clo$ and $g_n(y_n)$ tends to some point $y_{\infty} \in \clo$, where 
$x_{\infty} \neq y_{\infty}$. Let us fix a point $u \in \clo$ different from $a$. Up to a 
new subsequence, one can suppose that $g_n(u)$ tends to some point $b \in \clo$. We will 
prove that $(g_n)$ tends pointwise to this point $b \in \clo$ on $\clo \setminus \{ a \}$. 

Let us fix a point $v \in \clo$ different from $u$ and $a$. Let us first suppose that 
$v$ belongs to the interval $]u,a[$. In that case, for $n \in \mathbb{N}$ large enought 
one has $u\!<\!v\!<\!x_n\!<\!y_n\!<\!u$, and since $\nu$ is invariant
by $g \in G_{\nu}$, 
$$\nu([v,x_n] \times [y_n,u]) = \nu([g_n(v),g_n(x_n)] \times [g_n(y_n),g_n(u)]).$$
The left hand member of this equality converges to 
\hspace{0.05cm} $\nu([v,a[ \times ]a,u]) = \infty$
\hspace{0.03cm} as $n$ goes to infinity. Thus 
the right hand member also converges to infinity. But since $g_n(x_n)$, 
$g_n(y_n)$ and $g_n(u)$ converge to $x_{\infty}, y_{\infty}$ and $b$ respectively, and since 
$x_{\infty} \neq y_{\infty}$, this implies that $g_n(v)$ converges to $b$. When $v$ 
belongs to $]a,u[$, the same argument applied to the product $[u,x_n] \times [y_n,v]$ 
shows that $g_n(v)$ still converges to the point $b$.

A similar argument can be given for the sequence $(g_n^{-1})$, showing that $G_{\nu}$ has 
the convergence property. This finishes the proof of the proposition. 

\vspace{0.3cm}

The main example of a stable geodesic current is Liouville measure $Lv$ defined by 
$$dLv(u,v) = \frac{du \hspace{0.1cm} dv}{4 \sin^2(\frac{u-v}{2})}.$$
Other examples can be obtained by noting that stability is preserved by 
small perturbations. For example, for any function 
$K \in \mathcal{L}^2(\mathrm{S}^1 \times \mathrm{S}^1, Lv)$ which is 
invariant by the flip $(u,v) \mapsto (v,u)$, the Radon measure 
$\nu_K$ given by $d\nu_K = [1+K]^2 dLv$ is stable \cite{Na-kazh}. This 
kind of measures will be essential in what follows. 

%%%%%%%%%%%%%%%%%%%%%%%%%%%%%%%%%%%%%%%%%%%%%%%%%%%%%%%%%%%%%%%%%%%%%%%%%%%%%%%%%%%%%%%%%%%%%%%

\section{Reduction for Kazhdan's groups}

\subsection{Passing to the 3-fold covering} 

\hspace{0.45cm} The study of actions of Kazhdan groups on the circle is simplified by using a 
beautiful argument pointed to the author by D. Morris Witte: instead of considering the original 
action, it is sufficient (and easier) to deal with the induced action of a degree 3 central 
extension of the group on the 3-fold covering of the circle. However, it is not completely 
evident how to translate this argument to the setting of cocycles, and we will need to be 
a little bit careful to do that. The idea will consist on looking directly at an affine 
isometric action on some Hilbert space associated to this 3-fold covering. (Note that 
this argument does not appear in \cite{Na-superrig} or \cite{Na-kazh}, 
and it can be used to simplify some proofs therein.)

Let us denote by $\hat{\mathrm{S}}^1$ the 3-fold covering of the original circle $\clo$. For 
each $g \in G$ and a.e. $x \in \Omega$, the map $\alpha(g,x)\in \mathrm{Diff}_+^{1+\tau}(\clo)$ 
induces 3 diffeomorphisms $\hat{\alpha}_i(g,x)$, $i \in \{1,2,3\}$, of 
$\hat{\mathrm{S}}^1$. Those diffeomorphisms differ each one from the 
other by an order 3 Euclidean rotation. Remark that in general there is 
no canonical way to consider $\hat{\alpha}_i$ as a map from $G \times \Omega$ 
to $\mathrm{Diff}_{+}^{1+\tau}(\hat{\mathrm{S}}^1)$.

%%%%%%%%%%%%%%%%%%%%%%%%%%%%%%%%%%%%%%%%%%%%%%%%%%%%%%%%%%%%%%%%%%%%%%%%%%%%%%%%%%%%%%%%%%%%%%%%

\subsection{Construction of the affine isometric action} 

\hspace{0.45cm} Let us first consider the Hilbert space $\mathcal{H}'= 
\mathcal{L}^{2,\Delta}(\hat{\mathrm{S}}^1 \times \hat{\mathrm{S}}^1,
Lv)$ of square integrable real valued functions $K$ that satisfy a.e. the equalities 
$K(u,v)=K(v,u)$ and $K(u+2\pi/3,v+2\pi/3) = K(u,v)$. Now let us consider the 
Hilbert space $\mathcal{H}$ of functions $K \in \mathcal{L}^2(\Omega
\times \hat{\mathrm{S}}^1 \times \hat{\mathrm{S}}^1, \mu \times Lv )$
such that, for a.e. $x \in \Omega$, the function $K_x$ sending $(u,v)$
to $K(x,(u,v))$ belongs to $\mathcal{H}'$. 

Here is a main point. Even if there is not a well defined action of $G$ on $\Omega \times 
\hat{\mathrm{S}}^1$, one can naturally define the ``regular
representation" of $G$ on $\mathcal{H}$ by letting
$$\theta(g^{-1})K(x,(u,v)) = K(g(x),\hat{\alpha}_i(g,x)(u,v)) \hspace{0.06cm} 
[Jac(\hat{\alpha}_i)(g,x)(u,v)]^{1/2}.$$   
Indeed, this definition does not depend on $i \in \{1,2,3\}$. 
(Here and in what follows, we denote also by $\hat{\alpha}_i(g,x)$ the map from 
$\hat{\mathrm{S}}^1 \times \hat{\mathrm{S}}^1$ to itself induced diagonally by the 
original one of $\hat{\mathrm{S}}^1$, and we denote by $Jac(\hat{\alpha}_i)(g,x)(u,v)$ 
the Jacobian at the point $(u,v)$ of this map.) 

\vspace{0.25cm}

\begin{lem} {\em For each $g \in G$ the function $c(g)$ given by} 
$$(x,(u,v)) \longmapsto 1 - [Jac(\hat{\alpha}_i)(g^{-1},x)(u,v)]^{1/2}$$ 
{\em is well defined (i.e. it does not depend on $i \in \{1,2,3\}$) and 
belongs to the Hilbert space $\mathcal{H}$.}
\end{lem}

\noindent{\bf Proof.} The fact that $c(g)$ is well defined is clear. 
On the other hand, the proof of proposition 2.1 in \cite{Na-kazh}
shows that for each $g \in G$ and a.e. $x \in \Omega$, the map 
$$(u,v) \longmapsto 1 - [Jac(\hat{\alpha}_i)(g^{-1},x)(u,v)]^{1/2}$$ 
belongs to $\mathcal{H}'$. More precisely, there exists a constant 
$C_{\tau} < \infty$ such that
$$\left\| 1 - [Jac(\hat{\alpha}_i)(g^{-1},x)(u,v)]^{1/2} \right\| 
\leq C_{\tau} \thinspace \| \alpha_i'(g,x) \|_{\tau}.$$ 
Since by hypothesis the map $x \mapsto \alpha(g,x)$ takes essentially
its values in a bounded subset of $\mathrm{Diff}_+^{1+\tau}(\clo)$, 
this inequality proves the lemma. 

\vspace{0.3cm}

Now we have a candidate for an isometric affine action of $G$ on the Hilbert space 
$\mathcal{H}$, namely $A(g) = \theta(g) + c(g)$. The fact that this well defines 
an isometric action is straightforward to verify. 

%%%%%%%%%%%%%%%%%%%%%%%%%%%%%%%%%%%%%%%%%%%%%%%%%%%%%%%%%%%%%%%%%%%%%%%%%%%%%%%%%%%%%%%%%%%%%%%%

\subsection{The field of equivariant stable geodesic currents}

\hspace{0.45cm} When $G$ has Kazhdan's property (T), the isometric action $A$ above has 
a fixed point $K \in \mathcal{H}$, that is $\theta(g)K + c(g) = K$ for all $g \in G$. 
By definition, for all $g \in G$ one has a.e. the equality 
\begin{equation}
[1-K(x,(u,v))]^2=[1-K(g^{-1}(x),\hat{\alpha}_i(g^{-1},x)(u,v))]^2 
\hspace{0.08cm} Jac(\hat{\alpha}_i(g^{-1},x))(u,v).
\label{dos}
\end{equation}
Since $K$ belongs to $\mathcal{H}$, for a.e. $x \in \Omega$ the function 
$K_x$ that sends $(u,v)$ to $K(x,(u,v))$ belongs to $\mathcal{H}'$. We will denote by $\nu_x$ 
the measure on $\hat{\mathrm{S}}^1 \times \hat{\mathrm{S}}^1$ given by 
$$d\nu_x =[1-K_x]^2 dLv.$$ 
As a measure on $\hat{\mathrm{S}}^1 \times \hat{\mathrm{S}}^1 \setminus \Delta$, it is an 
a.e. well defined Radon measure. Moreover, by \S \ref{conv}, the measure $\nu_x$
is an a.e. well defined stable geodesic current. 

Let us denote by $\mathcal{SGC}(\hat{\mathrm{S}}^1)$ the space of Radon measures defined on 
$\hat{\mathrm{S}}^1 \times \hat{\mathrm{S}}^1 \setminus \Delta$ which
are stable, invariant by the flip $(u,v) \mapsto (v,u)$, and also 
invariant by the simultaneous order 3 Euclidean rotation 
$(u,v) \mapsto (u+2\pi/3,v+2\pi/3)$. 
The group $\mathrm{Homeo}_+(\clo)$ acts on $\mathcal{SGC}(\hat{\mathrm{S}}^1)$: 
for $h \in \mathrm{Homeo}_+(\clo)$ take one of its 3 preimages 
$\hat{h} \in \mathrm{Homeo}_+(\hat{\mathrm{S}}^1)$ and for 
$\nu \in \mathcal{SGC}(\hat{\mathrm{S}}^1)$ define $h(\nu)$ as 
$(\hat{h} \times \hat{h})_*(\nu)$. By the third property above, 
this definition is independent of the choice of the preimage $\hat{h}$. 

%%%%%%%%%%%%%%%%%%%%%%%%%%%%%%%%%%%%%%%%%%%%%%%%%%%%%%%%%%%%%%%%%%%%%%%%%%%%%%%%%%%%%%%%%%%%%%%%

\subsection{The reduction of the cocycle}

\hspace{0.45cm} We proved in \S 2.3 that there exists an equivariant map 
$\psi : \Omega \rightarrow \mathcal{SGC}(\hat{\mathrm{S}}^1)$, namely $\psi(x) = \nu_x$ 
(see equality (\ref{dos})). In \S \ref{liso} we will prove that the action of 
$\mathrm{Homeo}_+(\clo)$ on $\mathcal{SGC}(\hat{\mathrm{S}}^1)$ is smooth. Assuming 
this fact for a moment, 
we can apply the Cocycle Reduction Lemma to conclude that $\alpha$ is cohomologous in 
$\mathrm{Homeo}_+(\clo)$ to a cocycle $\beta : G \times \Omega \rightarrow H_{\nu}$, where $H_{\nu}$ 
is the stabiliser of some element $\nu \in \mathcal{SGC}(\hat{\mathrm{S}}^1)$ by the action of 
$\mathrm{Homeo}_+(\clo)$.\footnote{Is is maybe possible to perform this first reduction of the 
cocycle $\alpha$ inside the group of orientation preserving 
{\em quasisymmetric} homeomorphisms of the circle (compare with 
\cite{Hi2}, \cite{Na-superrig}, and the recent preprint \cite{mark}).}

Let us denote by $\hat{H}_{\nu}$ the degree 3 central extension of $H_{\nu}$. 
This is a subgroup of $\mathrm{Homeo}_+(\hat{\mathrm{S}}^1)$ which preserves $\nu$. So, by 
proposition (\ref{conver}), $\hat{H}_{\nu}$ is topologically conjugated to a subgroup of 
$\mathrm{PSL}(2,\mathbb{R})$. We claim that the action of $H_{\nu}$ on $\clo$ is free. 
Indeed, if $h \in H_{\nu}$ fixes one point of $\clo$, then one of its preimages $\hat{h}$ 
fixes 3 points of $\hat{\mathrm{S}}^1$. However, since $\hat{h}$ is topologically 
conjugated to an element of $\mathrm{PSL}(2,\mathbb{R})$, this is not possible 
unless $\hat{h}$ (and hence $h$) is the identity. 

We conclude that $H_{\nu}$ is a subgroup of $\mathrm{Homeo}_+(\clo)$ whose action is free and 
whose degree 3 central extension is topologically conjugated to a subgroup of the M\"oebius 
group. This is not possible unless $H_{\nu}$ is topologically conjugated to a group 
of Euclidean rotations. Modulo the proof of the smoothness of the action of $\mathrm{Homeo}_+(\clo)$ 
on $\mathcal{SGC}(\hat{\mathrm{S}}^1)$, this finishes the proof of Theorem A.

%%%%%%%%%%%%%%%%%%%%%%%%%%%%%%%%%%%%%%%%%%%%%%%%%%%%%%%%%%%%%%%%%%%%%%%%%%%%%%%%%%%%%%%%%%%%%%%%

\subsection{The smoothness of the action of $\mathrm{Homeo}_+(\clo)$ on 
$\mathcal{SGC}(\hat{\mathrm{S}}^1)$}
\label{liso}

\hspace{0.45cm} At first glance, this fact should seem rather surprising. For instance, it is 
easy to see that the action of $\mathrm{Homeo}_+(\clo)$ on the space of probability measures 
of the circle is not smooth. However, we will see that the stability property of the 
measures involved in our case are at the origin of many rigidity phenomena.

For the proof of the smoothness of the action of $\mathrm{Homeo}_+(\clo)$ on 
$\mathcal{SGC}(\hat{\mathrm{S}}^1)$, let us fix a sequence $(h_n)$ in $\mathrm{Homeo}_+(\clo)$ 
and an element $\nu\!\in\!\mathcal{SGC}(\hat{\mathrm{S}}^1)$ such that $\nu_n\!=\!h_n(\nu)$ 
tends to some $\nu_{\infty}\!\in\!\mathcal{SGC}(\hat{\mathrm{S}}^1)$. We have to prove that 
there exists $h \in \mathrm{Homeo}_+(\clo)$ such that $h(\nu) = \nu_{\infty}$. To do 
that, it suffices to prove that the sequences $(h_n)$ and $(h_n^{-1})$ are both 
equicontinuous. Indeed, in that case, any limit $h$ of a subsequence 
$(h_{n_k})$ of $(h_n)$ such that $(h^{-1}_{n_k})$ also converges will 
be an element of $\mathrm{Homeo}_+(\clo)$ satisfying $h(\nu)=\nu_{\infty}$.

Suppose by contradiction that $(h_n)$ is not equicontinuous and for each $n\in\mathbb{N}$ 
fix a preimage $\hat{h}_n \in \mathrm{Homeo}_+(\hat{\mathrm{S}}^1)$ of $h_n$. By passing 
to a subsequence if necessary, we may assume that there exist $\varepsilon > 0$ 
and two sequences of points $(u_n)$ and $(v_n)$ in $\hat{\mathrm{S}}^1$ such that 
$2\pi/3 > d(u_n,v_n) \geq \varepsilon$ for all $n \in \mathbb{N}$ and such that 
$d(\hat{h}_n^{-1}(u_n),\hat{h}_n^{-1}(v_n))$ converges to zero as $n$ goes to infinity. 
Again, by passing to a subsequence if necessary, we can suppose that $x_n=\hat{h}_n^{-1}(u_n)$ 
and $y_n = \hat{h}_n^{-1}(v_n)$ both converge to the same limit point $z\in\hat{\mathrm{S}}^1$, 
and that $u_n$ (resp. $v_n$) converges to some point $u_{\infty}$ (resp. $v_{\infty}$) in such a 
way that $2\pi /3 > d(u_{\infty},v_{\infty}) \geq \varepsilon$. Let us denote by $\bar{x}_n$, 
$\bar{\bar{x}}_n$, $\bar{y}_n$, $\bar{\bar{y}}_n$, etc, the points obtained from $x_n$, $y_n$, 
etc, by Euclidean rotations of order 3 (see Picture 1). 

\begin{figure}[hbtp]
\psfrag{a}{$x_n$} 
\psfrag{b}{$z$}
\psfrag{c}{$y_n$}
\psfrag{d}{$\bar{x}_n$}
\psfrag{e}{$\bar{z}$}
\psfrag{f}{$\bar{y}_n$}
\psfrag{g}{$\bar{\bar{x}}_n$}
\psfrag{h}{$\bar{\bar{z}}$}
\psfrag{i}{$\bar{\bar{y}}_n$}
\psfrag{j}{$\hat{h}_n$}
\psfrag{k}{$u_{\infty}$}
\psfrag{w}{$\hat{\mathrm{S}}^1$}
\psfrag{x}{$\hat{\mathrm{S}}^1$}
\psfrag{l}{$u_n$}
\psfrag{m}{$v_n$}
\psfrag{n}{$v_{\infty}$}
\psfrag{o}{$\bar{u}_{\infty}$}
\psfrag{p}{$\bar{u}_n$}
\psfrag{q}{$\bar{v}_n$}
\psfrag{r}{$\bar{v}_{\infty}$}
\psfrag{s}{$\bar{\bar{u}}_{\infty}$}
\psfrag{t}{$\bar{\bar{u}}_n$}
\psfrag{u}{$\bar{\bar{v}}_n$}
\psfrag{v}{$\bar{\bar{v}}_{\infty}$}
\centerline{\includegraphics{coc.eps}}
\label{coc}
\end{figure}

For all $n \in \mathbb{N}$ one has 
\begin{equation}
\nu ([y_n,\bar{x}_n] \times [\bar{y}_n,x_n]) = \nu_n ([v_n,\bar{u}_n] \times [\bar{v}_n,u_n]).
\label{ig}
\end{equation}
For each $k \in \mathbb{N}$ let us fix 4 points $p_k,q_k,\bar{p}_k,\bar{q}_k$ such that 
$p_k\!<z\!<q_k\!<\!\bar{p}_k\!<\!\bar{z}\!<\!\bar{q}_k$ and such that 
$d(p_k,z)=d(z,q_k)=d(\bar{p}_k,\bar{z})=d(\bar{z},\bar{q}_k)=1/k$. For fixed $k$ there exists 
a positive integer $n(k)$ such that for all $n \geq n(k)$ one has $d(x_n,z)<1/k$ and $d(z,y_n)<1/k$. 
So, for $n \geq n(k)$, 
$$\nu([y_n,\bar{x}_n]\times [\bar{y}_n,x_n]) \geq \nu([q_k,\bar{p}_k]\times [\bar{q}_k,p_k]).$$ 
Thus
$$\liminf_{n \rightarrow \infty} \nu([y_n,\bar{x}_n]\times [\bar{y}_n,x_n]) \geq 
\nu([q_k,\bar{p}_k]\times [\bar{q}_k,p_k]).$$
Since this is true for all $k \in \mathbb{N}$, one obtains 
$$\liminf_{n \rightarrow \infty} \nu([y_n,\bar{x}_n]\times [\bar{y}_n,x_n]) \geq 
\nu(]z,\bar{z}[ \times ]\bar{z},z[) = \infty.$$
However, the right hand term of equality (\ref{ig}) tends to 
$\nu_{\infty}([v_{\infty},\bar{u}_{\infty}] \times [\bar{v}_{\infty},u_{\infty}])$, 
and the value of this expression is finite. (Recall that $\nu_{\infty}$ is a Radon 
measure on $\hat{\mathrm{S}}^1\times \hat{\mathrm{S}}^1 \setminus \Delta$.) This 
contradiction finishes the proof of the equicontinuity of $(h_n)$. A similar argument 
allows to show the equicontinuity of $(h_n^{-1})$, concluding the proof of the smoothness 
of the action of $\mathrm{Homeo}_+(\clo)$ on $\mathcal{SGC}(\hat{\mathrm{S}}^1)$.

%%%%%%%%%%%%%%%%%%%%%%%%%%%%%%%%%%%%%%%%%%%%%%%%%%%%%%%%%%%%%%%%%%%%%%%%%%%%%%%%%%%%%%%%%%%%%%%%

\section{The higher rank case}

\hspace{0.45cm} The proof of Theorem B is obtained by putting the affine isometric action $A$ 
of previous paragraphs into the context of \S 1.3. The case in which this action has a global 
fixed point was ruled out in \S 2. If this is not the case but the action is non uniform, 
then we are in the ``degenerated case", and we can apply Delorme's Lemma. Finally, in 
the uniform case, we will use Shalom's Superrigidity Theorem.

%%%%%%%%%%%%%%%%%%%%%%%%%%%%%%%%%%%%%%%%%%%%%%%%%%%%%%%%%%%%%%%%%%%%%%%%%%%%%%%%%%%%%%%%%%%%%%%%

\subsection{The degenerated case} 

\hspace{0.45cm} Associated to the cocycle 
$\alpha: G \times \Omega \rightarrow \mathrm{Diff}_+^{1+\tau}(\clo)$, 
$\tau > 1/2$, let us consider the affine isometric action $A$ on the Hilbert space 
$\mathcal{H}$. The arguments in \S 2 show that if this action has an 
invariant vector $K \in \mathcal{H}$, then $\alpha$ is cohomologous in 
$\mathrm{Homeo}_+(\clo)$ to a cocycle taking values in the group of Euclidean rotations.

Let us suppose now that there is no invariant vector for $A$ but this action is not uniform. 
By Delorme's Lemma, there exists a sequence $(K_n)$ of unitary vectors in $\mathcal{H}$ 
such that for all $g \in G$ one has \hspace{0.03cm} 
$\lim_{n \rightarrow +\infty} \| \theta(g)K_n - K_n \| = 0$. 
For each $n \in \mathbb{N}$ let us consider the probability 
measure $\hat{m}_n$ on $\Omega \times \hat{\mathrm{S}}^1$ defined by  
$$\frac{d \hat{m}_n}{d(\mu \times Leb)}(x,u) = \int_{\hat{\mathrm{S}}^1} \frac{K_n^2(x,(u,v))}
{4 \sin^2(\frac{u-v}{2})} \hspace{0.08cm} dv.$$
Up to a subsequence, one can suppose that $\hat{m}_n$ tends to some limit $\hat{m}$. 
This measure $\hat{m}$ induces a probability measure $m$ on $\Omega \times \clo$ by projection 
$\Omega \times \hat{\mathrm{S}}^1 \rightarrow \Omega \times \clo$, and a 
straightforward computation shows that $m$ is invariant by the skew action of G.

%%%%%%%%%%%%%%%%%%%%%%%%%%%%%%%%%%%%%%%%%%%%%%%%%%%%%%%%%%%%%%%%%%%%%%%%%%%%%%%%%%%%%%%%%%%%%%%%%

\subsection{The uniform case}

\hspace{0.45cm} If the affine action $A$ is uniform, then Shalom's Superrigidity Theorem gives 
the existence of some index $i \in \{1,\ldots,k\}$ such that the space $\mathcal{H}_i$ of 
$\theta(G_i')$-invariant vectors is non trivial. Due to the commutativity of $G_i$ and 
$G_i'$, this space $\mathcal{H}_i$ is $\theta(G)$-invariant.

Each non zero function $K \in \mathcal{H}_i$ gives rise to an almost everywhere defined 
family of finite measures on the fibers over each point $x \in \Omega$. Indeed, 
denoting by $\pi : \hat{\mathrm{S}}^1 \rightarrow \clo$ the canonical 
projection, one can define, for $X \subset \clo$, 
$$\mu_{K,x}(X) = \int_{\hat{\mathrm{S}}^1} \int_{\pi^{-1}(X)} 
\frac{K^2(x,(u,v))}{4\sin^2(\frac{u-v}{2})} \hspace{0.07cm} du \hspace{0.05cm} dv.$$ 

Let us denote by $\mathcal{F}_i$ the collection of all families of finite measures obtained from 
$\mathcal{H}_i$ by this procedure, and normalised so that 
$$\int_{\Omega} \mu_{K,x}(\clo) \hspace{0.05cm} d\mu(x)=1.$$ 
If $\mathcal{F}_i$ contains only one element, then we obtain a normalised equivariant 
family of a.e. finite measures on $\clo$, and so an invariant probability measure for the skew 
action of $G$ on $\Omega \times \clo$. Suppose in what follows that $\mathcal{F}_i$ contains 
more than one element. 

Let us fix an orthonormal basis $\{K_1,K_2,\ldots\}$ of $\mathcal{H}_i$, and let us define 
$$\bar{K} = \sum_{i} \frac{|K_i|}{2^i}, \qquad K = \frac{\bar{K}}{\| \bar{K} \|}.$$
This function $K \in \mathcal{H}_i$ gives an element of $\mathcal{F}_i$ such that the 
support of $\mu_{K,x}$ is ``maximal" on the fiber of a.e. $x \in \Omega$. Denote by 
$\mathrm{S}^1_K$ the equivariant topological circle obtained by cutting up the connected 
components of the complementary of the support of measures $\mu_{x,K}$. (Note that 
equivariance follows from the fact that $\mathcal{H}_i$ is a $\theta(G)$-invariant 
subspace.) Denote also by $\mu_K$ the ($G_i'$-equivariant) probability measure 
on the ($G$-equivariant topological) circle $\clo_K$ induced by the family 
$\mu_{x,K}$, with $x \in \Omega$. The procedure of cutting induces a new cocycle 
$$\bar{\alpha}: G \times \Omega \rightarrow \mathrm{Homeo}_+(\clo_K),$$ 
which can be considered as a ``semiconjugated" of the cocycle $\alpha$. 

By hypothesis, there exists a unitary $K' \in \mathcal{H}_i$ inducing a different 
element of $\mathcal{F}_i$ from that induced by $K$. Denote by $\mathcal{FM}(\clo_K)$ the 
space of finite measures on $\clo_K$, and denote by $\mathcal{LM}(\clo_K)$ the subspace 
of $\mathcal{FM}(\clo_K)$ whose elements have no atoms and total support ({\em i.e.} the 
measure of each non empty open set is positive). Let $\mathcal{M}^{\Delta}(\clo_K)$ be 
the space $\mathcal{LM}(\clo_K)\!\times\!\mathcal{FM}(\clo_K) \setminus \Delta$, where 
$\Delta$ denotes the corresponding diagonal. The map $\psi: G_i' \rightarrow 
\mathcal{M}^{\Delta}(\clo_K)$ given by $\psi(g_i') = (\mu_{K},\mu_{K'})$ 
is $G_i'$-equivariant, and it is easy to verify that the action of 
$\mathrm{Homeo}_+(\clo_K)$ on $\mathcal{M}^{\Delta}(\clo_K)$ is smooth. 
So, applying the Cocycle Reduction Lemma to the restriction $\bar{\alpha}|_{G_i'}$, 
we obtain that $\bar{\alpha}$ is cohomologous to some cocycle 
$\bar{\beta} : G \times \Omega \rightarrow \mathrm{Homeo}_+(\clo_K)$ 
which takes values in the stabiliser of some point $(\mu,\mu')$ in 
$\mathcal{LM}(\clo_K)\!\times\!\mathcal{FM}(\clo_K) \setminus \Delta$. On the other 
hand, the intersection $H$ of the stabilisers 
of $\mu$ and $\mu'$ is a finite subgroup of $\mathrm{Homeo}_+(\clo)$. 
Indeed, the stabiliser of $\mu$ is topologically conjugated to the
group of Euclidean rotations, and the only closed and strict subgroups of the group of 
rotations are the finite ones.

By identifying (equivariantly) the points of the orbits by $H$, one obtains another 
topological equivariant circle $\clo_K /\!\sim$, and $\bar{\beta}$ induces a cocycle 
$$\beta : G \times \Omega \rightarrow \mathrm{Homeo}_+(\clo_K /\!\sim)$$ 
such that the image $\beta(G_i' \times \Omega)$ is trivial. (Thus $\bar{\beta}$ 
can be seen as a ``finite cover" of the cocycle $\beta$.) Since $G_i$ and $G_i'$ 
commute, for $g \in G_i$, $g' \in G_i'$ and a.e. $x \in \Omega$, we have  
$$\beta(g,g'(x)) = \beta(g,g'(x))\beta(g',x) = \beta(gg',x) = \beta(g'g,x) = 
\beta(g',g(x)) \beta(g,x) = \beta(g,x).$$
Since $G_i'$ acts ergodically on $\Omega$, we conclude that $\beta$
does not depend on the point $x \in \Omega$, and thus it is indeed a Borel 
(hence continuous) homomorphism from $G$ to $\mathrm{Homeo}_+(\clo_K/\!\sim)$ 
factoring through $G_i$. This finishes the proof of the first part of Theorem B. 

Suppose now that each $G_i$ is non discrete and almost topologically simple ({\em i.e.} the  
closed normal subgroups are compact or cocompact). Recall that each locally compact 
subgroup of $\mathrm{Homeo}_+(\clo)$ is a (real) Lie group \cite{ghys2}. Moreover, if 
such a group is connected, then it is topologically conjugated to a subgroup of a direct 
product of some groups of the following list: $\mathbb{R}$, $\mathrm{Aff}_+(\mathbb{R})$, 
$\mathrm{SO}(2,\mathbb{R})$, $\mathrm{PSL}_k(2,\mathbb{R})$ for some $k \geq 1$ or 
$\widetilde{\mathrm{PSL}}(2,\mathbb{R})$. The cocycle $\beta$ then induces naturally  
homomorphisms from $G_i$ to some of these groups, and $G_i$ being non discrete and almost 
topologically simple, those homomorphisms have compact image or are surjective with compact 
kernel. It is easy to see that the only possibility for which there is no invariant 
probability measure for the $G$-action on $\Omega \times \clo$ is that $\beta(G_i)$ 
is topologically conjugated to some finite cover 
$\mathrm{PSL}_k(2,\mathbb{R})$ of $\mathrm{PSL}(2,\mathbb{R})$, 
and this finishes the proof of Theorem B.

\vspace{0.1cm}

\begin{rem} When each $G_i$ is a simple Lie group, it has been proven in \cite{WZ} 
that the invariant probability measure on $\Omega \times \clo$ can be taken so 
that its projection into $\Omega$ coincides with the original measure $\mu$. 
We ignore if this is still true in our setting.
\end{rem}

%%%%%%%%%%%%%%%%%%%%%%%%%%%%%%%%%%%%%%%%%%%%%%%%%%%%%%%%%%%%%%%%%%%%%%%%%%%%%%%%%%%%%%%%%%%

\section{Appendix: improving the differentiability class}

\hspace{0.45cm} Recently, some new rigidity phenomena have been discovered for Kazhdan's 
groups. More precisely, the fixed point property is not only valid for isometric actions 
of such groups on Hilbert spaces, but also on any $L^p$-space for $p \in [1,2]$ and for 
``almost" isometric actions: see \cite{BG} and \cite{FM} respectively. By combining the 
methods of those works, one can easily prove that for any group $G$ satisfying Kazhdan's 
property (T), there exists a constant $\delta(G) > 0$ such that $G$ satisfies the fixed 
point property for isometric actions on $L^p$-spaces for all $p \in [1,2+\delta(G)]$.

By considering an $L^p$-version of Liouville's cocycle and then using the same arguments of 
\cite{Na-kazh} and this work, the preceding result allows to prove a sharper version of 
theorem A for cocycles into $\mathrm{Diff}_+^{1+\tau}(\clo)$ for $\tau > \tau(G)$, 
where $\tau(G) < 1/2$ depends on $G$ (one can take 
$\tau(G) = 1/\delta(G)$).\footnote {It is unknown 
but very plausible that Theorem A (or at least the 
main result of \cite{Na-kazh}) still holds for $C^1$-diffeomorphisms. A very 
interesting problem is to study the case of piecewise affine homeomorphisms 
of the circle. For instance, it is not known if some generalized Thompson's 
groups satisfy Kazhdan's property (T) (see however \cite{farley}).} A similar 
(and stronger) remark can be done for the higher rank case. Indeed, it has 
been recently remarked \cite{BG} that a weak version of the Superrigidity 
Theorem for Reduced Cohomology holds for representations 
in $L^p$-spaces for any $p \in ]1,\infty[$. Using this new 
result, it is not very difficult to obtain an improved version of Theorem 
B (and the related results from \cite{Na-superrig}) for cocycles (resp. for 
representations) into $\mathrm{Diff}_+^{1+\tau}(\clo)$ for any $\tau > 0$.

\vspace{0.3cm}

\noindent{\bf Acknowledgments.} The author would like to thank A. Valette for his 
invitation to the University of Neuch\^atel (where this work was carried out), and 
D. Gaboriau for many interesting discussions on the subject.  

\vspace{0.1cm}

%%%%%%%%%%%%%%%%%%%%%%%%%%%%%%%%%%%%%%%%%%%%%%%%%%%%%%%%%%%%%%%%%%%%%%%%%%%%%%%%%%%%%%%%%%%%

\begin{small}

\vspace{0.2cm}

\noindent Andr\'es Navas\\

\vspace{0.1cm}

\noindent IH\'ES, 35 route de Chartres, 91440 Bures sur Yvette, France 
(anavas@ihes.fr)\\

\vspace{0.1cm}

\noindent Univ. de Chile, Las Palmeras 3425, $\tilde{\mathrm{N}}$u$\tilde{\mathrm{n}}$oa, 
Santiago, Chile (andnavas@uchile.cl)\\

\end{small}

\end{document}